\documentclass{svjour3}

\usepackage{amscd,amsfonts,amssymb,amsmath,latexsym,array,hhline,xcolor,graphicx}

\usepackage{latexsym}
\usepackage{times}

\newcommand\F{\mbox{I\kern-2pt F}}

\newcommand\cE{{\cal E}}

\newcommand\cF{{\cal F}}

\newcommand\cL{{\cal L}}

\newcommand\ups{{\upsilon}}
\newcommand\Ups{{\Upsilon}}
\newcommand\e{{\varepsilon}}
\newcommand{\wt}{\widetilde}

\def\E{{\bf E}}

\def\P{{\bf P}}

\def\R{{\mathbb R}}
\def\bbr{{\mathbb R}}
\def\N{{\mathbb N}}

\newcommand\fdem{$\Box$}
\newcommand\beq{\begin{equation}}
\newcommand\eeq{\end{equation}}
\newcommand\bea{\begin{eqnarray}}
\newcommand\eea{\end{eqnarray}}
\newcommand\bean{\begin{eqnarray*}}
\newcommand\eean{\end{eqnarray*}}

\begin{document}
\title{On  ruin probabilities with investments  in a risky asset  with a switching regime price}

\author{Yuri Kabanov  
\and Serguei\ ~ Pergamenshchikov}

\date{}
\institute{
  \at
              Laboratoire de Math\'ematiques, Universit\'e Bourgogne 
Franche-Comt\'e,
Besan\c{c}on, France, and 
Lomonosov Moscow State University, Moscow, Russia\\
 \email{Youri.Kabanov@univ-fcomte.fr}    \\
\and Laboratoire de Math\'ematiques Rapha\"el Salem, 
Universit\'e de Rouen, 
France,\\  and 
National Research Tomsk State University,  Russia \\
  \email{Serge.Pergamenchtchikov@univ-rouen.fr}   }

\date{Received: date / Accepted: date}

\titlerunning{On  ruin probabilities  }

\maketitle

\begin{abstract}
We investigate  the asymptotic of ruin probabilities when the company    invests its reserve in a risky asset  with a switching regime price. We assume that the asset price is a conditional geometric Brownian motion with parameters modulated by a Markov process with a finite number of states. 
 Using the technique of the implicit renewal theory 
we obtain the rate of convergence to zero of the ruin probabilities as the initial capital tends to infinity.  
\end{abstract}



 
   \section{Introduction}  
 Models, where an insurance company 
invests its reserve  (or a part of it) in a risky asset, constitutes an important class currently under an extensive study. Considering  a single risky asset is justified by the common practice of investing in a market portfolio 
or in an index (a fund which simulates an index like the DAX or the S\&P500) which is an economically reasonable strategy. Since the insurance contacts usually are of a long duration  and the return may depend, e.g., on the  business cycles of economy, models with regime switching  are now more and more popular.  The main question is the rate of decay of the ruin probability as the initial reserve tends to infinity. 

In this note 
we extend the recent result of Ellanskaya and Kabanov, \cite{EK}, established 
for a model with   characteristics of the asset price depending on  a telegraph process, i.e. on a Markov process with two states, $0$ and $1$.   Here we study  the case where the characteristics depend on the ergodic Markov process
$\theta$ with a finite 
number  of states.  When $\theta_t=k$, the asset price evolves as a  geometric Brownian motion with drift  $a_k$ and volatility
$\sigma_k$. It is well known, see, e.g. \cite{FrKP,KP,KPukh}, that in the case of a single regime, i.e. when the price process is the classic  gBm with drift  $a$ and   volatility $\sigma$,
the ruin probability decreases to zero as the initial capital $u$ tends to infinity with the rate $\beta:=2a/\sigma -1$  provided that $\beta>0$. In \cite{EK} is was shown that in the model with two regimes,  $0$ and $1$, the ruin 
probabilities decrease with a rate $\beta$ where $\beta$ is a number between the  values  $\beta_k:=2a_k/\sigma_k -1$, $k=0,1$, assumed to be strictly positive. This $\beta$ is  the root of an algebraic   equation of third order and does not depend on the initial value of $\theta$. 

In the present paper we extend this result to the case where
the number of states of the hidden Markov process $\theta$ is $K\ge 2$. It happens that,  provided all $\beta_k>0$, $k=0,\dots, K-1$, the rate of convergence
to zero of the ruin probabilities, depending, in general,  on the initial state $i$, is a root of the cumulant generating function of  the value of   log price process at the first return time of $\theta$ to the state $i$. 
It is worth to note that  the switching by  telegraph signal is rather specific: the latter returns to the initial state after the second jump while for  a general Markov process the return may happen after arbitrary number of jumps. 

Though the main idea is again based on the implicit renewal theory, it happens that the analysis of the considered  model  is much more complicated and the calculation of the rate parameter, depending, in 
general, on the initial state, is not so straightforward. 
We hope that the result of this paper elucidate challenging problems of estimation of ruin probabilities for other stochastic volatility models.  


\section{The model}
Let   $(\Omega,\cF,{\bf F}=(\cF_t)_{t\in \R_+},\P)$  be a stochastic basis with a Wiener process $W=(W_t)$, 
a  Poisson random measure $\pi(dt,dx)$ on $\R_+\times \R$ with 
the mean $\tilde \pi(dt,dx)=\Pi(dx)dt$, and a piecewise constant right-continuous  Markov process $\theta=(\theta_t)$. For the latter we assume that it  takes values in the finite set  $\{0,1,...,K-1\}$,  has  the $K\times K$
transition intensity matrix $\Lambda=(\lambda_{jk})$ with communicating states,   
and the initial value $\theta_0=i$ (so that $\theta=\theta^i)$.     
The $\sigma$-algebras generated by $W$, $\pi$, and $\theta$ are  independent. 

Recall that $\lambda_{jj}=-\sum_{k\neq j}\lambda_{jk}$ and $\lambda_i:=-\lambda_{ii}>0$ for each $i$. 

 Let  $T_n$ be the successive jumps  
of the Poisson process $N=(N_t)$ with $N_t:=\pi([0,t],\R)$ and let  $\tau_n$ be  the successive jumps  of $\theta$ with the convention  $T_0=0$ and $\tau_0=0$. 

Recall that the lengths of the intervals between the consecutive jumps of $\theta$ are independent  exponentially distributed random variables.

The  reserve 
$X=X^u$ of an insurance company evolves not only due to the business activity 
part, described as in the classical 
Cram\'er--Lundberg model,  but also due to the stochastic interest rate.  We assume that  the reserve is fully invested  in a risky asset whose price 
$S$ is  a conditional geometric Brownian motion given the Markov process $\theta$. That is,  $S$ is given by a so-called hidden Markov model with 
$$
dS_t=S_t (a_{\theta_t} dt+\sigma_{\theta_t}dW_t), \qquad S_0=1,   
$$ 
where $a_k\in \R$, $\sigma_k>0$, $k=0,...,K-1$.   In this case,
$X$ is of the form
\beq
 \label{risk}
 X_t=u+  \int_0^t  X_s  dR_s+ dP_t   
 \eeq
 where $dR_t=a_{\theta_t}dt+\sigma_{\theta_t} dW_t=dS_t/S_t$, that is  $R$ 
 is the relative price  process,   
 and 
\beq
 \label{r}
P_t=ct  + \int_0^t\int x  \pi(dt,dx)=ct+x*\pi_t. 
\eeq

So, the  reserve evolution  is described by the process $(X^{u},\theta)=
(X^{u,i},\theta^i)$ 
where 
  $u>0$ is the initial capital and $i$ is the initial regime, i.e. the initial value of $\theta$. 

We assume that  $P$ is not an increasing process: otherwise the probability of ruin is zero. 

We also assume that 
$\Pi(\R)<\infty$, that is 
$\Pi(dx)=\alpha_1 F_1(dx)+\alpha_2F_2(dx)$ where 
$F_1(dx)$ is a probability distribution on $]-\infty,0[$ and $F_2(dx)$ is a probability distribution on $]0,\infty[$. In this
case the integral with respect to the jump measure is simply a difference of two independent  
compound Poisson processes with intensities $\alpha_1$, $\alpha_2$ of  jumps  downwards and upwards and whose absolute values have  the distributions $F_1(dx)$ and $F_2(dx)$, respectively. 

The solution of the linear equation (\ref{risk}) can be represented as   
\beq
\label{uY}
X_t^{u}=\cE_t(R)(u-Y_t)=e^{V_t}(u-Y_t) 
\eeq 
where  
\beq
\label{Y_t}
Y_t:=-\int_{[0,t]} \cE^{-1}_{s}(R)dP_s=-\int_{[0,t]}  e^{-V_{s}}dP_s=-e^{-V}\cdot P_t, 
\eeq
the stochastic exponential $\cE_t(R)$ is equal to   ${S_t}$, and  
the  log price process $V=\ln \cE(R)$ admits the stochastic differential  
 $$
 dV_s=\sigma_{\theta_s}dW_s+(a_{\theta_s}-(1/2)\sigma_{\theta_s}^2)ds, \quad V_0=0. 
 $$  
 Of course,  $S$, $R$, $Y$, and $V$ depend on  $i$ (we omitted the superscript $i$ in the above formulae).   

\smallskip 
 Let $\tau^{u,i}:=\inf \{t>0:\ X^{u,i}_t\le 0\}$ be the instant of ruin corresponding to the initial capital $u$ and the initial regime $i$. 
 Then 
 $\Psi_i(u):=\P[\tau^{u,i}<\infty]$ is the ruin probability  
 and $\Phi_i(u):=1-\Psi_i(u)$
 is the survival probability. 
It is clear that   $ \tau^{u,i}= \inf \{t\ge 0:\ Y^ i_{t} \ge u\}$.

Recall that the constant parameter values $a=0$, $\sigma=0$, correspond to the 
Cram\'er--Lundberg  model for which  the process $X^{u,i}=u+P_t$.  In the actuarial literature the compound Poisson process $P$   is usually written  in the form 
\beq
\label{r0}
P_t= ct  - \sum_{k=1}^{N_t}\xi_k
\eeq
where  either $\xi_k\ge 0$, $c> 0$ (i.e. $F_2=0$ --- jumps only downwards --- the case of non-life insurance)  or $\xi_k\le 0$, $c< 0$ (i.e. $F_1=0$ --- jumps only upwards --- the case of life insurance or annuity payments).  Models with both kinds of jumps 
are also  considered in the literature, see, e.g., \cite{AGY} and references therein.  
 For the classical models  with a positive average trend and $F$
having a ``non-heavy" tail, the Lundberg inequality asserts that the ruin probability decreases
exponentially as the initial capital  $u$ increases  to infinity.
For the exponentially distributed claims the ruin
probability admits an explicit expression, see \cite{Asm},  Ch. IV.3b, or
\cite{Gr}, Section 1.1. 

For the models with  investment in a risky asset the situation is completely 
different. For example, for the model with exponentially distributed jumps 
and  price following a geometric Brownian motion with the drift coefficient 
$a$ and the volatility $\sigma>0$ in the case where $2a/\sigma^2-1>0$, the ruin probability as a function of the initial capital $u$, decreases as $C u^{1-2a/\sigma^2}$. If  
$2a/\sigma^2-1\le 0$ the ruin happens with probability one, see \cite{FrKP}, \cite{KP}, \cite{PZ}, \cite{KPukh}.

To formulate our result for the model where the volatility and drift are  modulated by 
a finite-state Markov process {\bf we assume} throughout this note that 
except  Section \ref{sec:prob1} that 
\beq
\label{beta}
\beta_j:=2a_j/\sigma_j^2-1>0, \quad j=0,...,K-1,
\eeq
(in other words, $\beta_*:=\min_j \beta_j>0$). 
\smallskip

Let $\ups^ i_1:=\inf \{t>0\colon \theta^{i}_{t-}\neq i, \ \theta^{i}_t=i\}$ be the first return time of the (continuous-time) Markov process $\theta=\theta^i$ to its initial state $i$. We consider further  the consequent return times defined recursively: 
$$
\ups^ i_k:=\inf \{t>\ups^i_{k-1}\colon \theta^{i}_{t-}\neq i, \ \theta^{i}_{t}=i\}, \quad k=2,... 
$$

We introduce the  random variable $M_{i1}:=e^{-V_{\ups^i_1}}$ and define   the 
moment  generating function $\Ups_i :\R_+\to \bar \R_+:=\R_+\cup \{\infty\}$ with 
$$
\Ups_i(q):= \E [M^q_{i1}].
$$
\begin{proposition} 
\label{prop1}
The function $\Ups_i$ is strictly convex, continuous, and there is unique $\gamma_i>0$ such that  $\Ups_i(\gamma_i)=1$.  
\end{proposition}

Note that one can characterize $\gamma_i$ also as the strictly positive roots of the cumulant  generating functions $H_i(q):=\ln \E e^{-V_{\ups^i_1}}$, strictly convex and continuous.    
 
Postponing the proof of Proposition \ref{prop1} to the next section we formulate our main result: 

\begin{theorem}  
\label{main}
Suppose that $\Pi(|x|^{\gamma_i}):=\int |x|^{\gamma_i} \Pi(dx) <\infty$. Then 
$$
0<\liminf_{u\to \infty} u^{\gamma_i} \Psi_i(u)\le \limsup_{u\to \infty}u^{\gamma_i} \Psi_i(u)<\infty. 
$$
  \end{theorem}

{\sl Important remark.} In  the case where $\theta$ is a telegraph signal, i.e. 
a two-state Markov process,  the values $\gamma_0$ and $\gamma_1$ coincide (see \cite{EK}).  In the considered general case, $\gamma_i$ may depend on the initial value $i$. 
To alleviate formulae, we fix the initial value $i=0$  and  omit  the index $i=0$ when this  does not lead to ambiguity.  
\smallskip 

The proof of Theorem \ref{main}  is based on the implicit renewal theory.   
 To apply  it, we verify that the random variables $
Q=Q_1:=- e^{-V}\cdot P_{\ups_1} $  
belong to $L^{\gamma}(\Omega)$, the process  $Y$ has at infinity a finite limit $Y_\infty$ and show that $Y_\infty$ is a random variable unbounded from above with  the same law as $Q+MY_{\infty}$ where 
$M=M_{1}$ and  $\E [M_{1}^{\gamma}]=1$.  Also  $\E[ M_1^{\gamma+\e}]<\infty$ for some 
$\e>0$.  Clearly, the law of the random variable $\ln M_{1}$ is not arithmetic. We establish the bounds  $\bar G(u)\le \Psi_i (u)\le  C \bar G(u)$ where $\bar G_i(u)=\P[Y_\infty >u]$ and a constant $C>0$, Lemma \ref{G-Paulsen}. 

With these  facts Theorem \ref{main}  follows from Theorem \ref{Le.3.2}  below which is the Kesten--Goldie theorem, see Th. 4.1 in \cite{Go91}, combined with a statement on strict positivity of $C_+$ due to Guivarc'h and Le Page, \cite {BD} (for a simpler proof of the latter  see    Buraczewski and Damek,  \cite{BD}, and an extended discussion in Kabanov and Pergamenshchikov, \cite{KP2020}). \begin{theorem} 
\label{Le.3.2} 
Let $Y_{\infty}$ has the same law as $Q+MY_{\infty}$ where $M>0$.  
Suppose that $(M,Q)$ is such that the law of   $\ln\,M$ 
is non-arithmetic and, for  some $\gamma>0$,
\begin{align}\label{3.3}
\E [M^\gamma]=1, \ \ \ \E[ M^\gamma\,(\ln\,M)^+]<\infty, \  \  \ 
\E[ |Q|^\gamma ]<\infty. 
\end{align}
Then 
\bean
&&\lim_{u\to\infty}\,u^\gamma\,\P[Y_{\infty} >u]=C_+<\infty,\\
&&\lim_{u\to\infty}\,u^\gamma\,\P[Y_{\infty} <-u)]=C_-<\infty,
\eean
where $C_++C_->0$. 

If the  random variable  $Y_{\infty} $ is unbounded from above, then $C_+>0$. 
\end{theorem} 

\section{Properties of the moment generating function: the  proof of Proposition \ref{prop1}}  
\label{sec3}
Recall that $\tau_n$ are the  moments of consecutive jumps of $\theta$, that is, $\tau_0:=0$, 
$$
\tau_n:=\inf\{t>\tau_{n-1}\colon \theta_{t-}\neq \theta_t\}, \quad n\ge 1.
$$   

We introduce  the imbedded Markov chain $\vartheta_n:=\theta_{\tau_n}$, $n=0,1,...$ with 
transition probabilities $P_{kl}=\lambda_{kl}/\lambda_{k}$, $k\neq l$, and $P_{kk}=0$.   
Then  
$\varpi:=\inf\{j\ge 2\colon \vartheta_j=0\}$ is the first return time of the (discrete-time) Markov chain  $\vartheta$ to the starting point $0$ and  $\ups_1=\tau_\varpi$. 

Put 
$$
Z^j_t:=\sigma_j W_t+(a_j-\sigma_j^2/2) t=\sigma_j W_t+(1/2)\sigma_j^2\beta_j t. 
$$

The random variable $M_{1}$ admits the representation 
$$
 M_{1}=\sum_{k\ge 2} \sum_{i_1\neq 0,i_2\neq 0,...,i_k=0}
 I_{\{\vartheta_1=i_1,\vartheta_2=i_2,...,\vartheta_k=i_k\}}e^{-\zeta^{0}_1}e^{-\zeta^{i_1}_2}...e^{-\zeta^{i_{k-1}}_k}
$$ 
where 
$$
\zeta^{0}_{1}:=Z^{0}_{\tau_1}-Z^{0}_{\tau_0},\quad  \zeta^{i_1}_{2}:=Z^{i_1}_{\tau_2}-Z^{i_1}_{\tau_1}, \quad ....,\quad  \zeta^{i_{k-1}}_{k}:=Z^{i_{k-1}}_{\tau_k}-Z^{i_{k-1}}_{\tau_{k-1}}.
$$
The conditional law of   random variables  $\zeta^{0}_1$,...,$\zeta^{i_{k-1}}_k$ given $\vartheta_1=i_1,\vartheta_2=i_2,...,\vartheta_k=i_k$
is the same as the unconditional law of independent random variables $\tilde \zeta^{0}_1$,...,$\tilde \zeta^{i_{k-1}}_k$.  For any $m$ the law $\cL(\tilde \zeta^{j}_m)=\cL(\sigma_j W_\tau +(1/2)\sigma_j^2\beta_j\tau)$ where 
an exponential random variable $\tau$ with  parameter $\lambda_{j}$ is independent on the Wiener process $W$. 

It follows that 
\beq
\label{H0}
\Ups(q):=\E [ M^q_{1}]=\sum_{k\ge 2} \sum_{i_1\neq 0,i_2\neq 0,...,i_k=0}P_{0i_1}P_{i_1i_2}...
P_{i_{k-1}0}f_0(q)f_{i_1}(q)...f_{i_{k-1}}(q),
\eeq
 where 
 $$
f_{j}(q)= \E[ e^{-q\tilde \zeta^{j}_m}]=\lambda_{j}
\E\left[  \int_{0}^\infty e^{-q(\sigma_j W_t+(1/2)\sigma_j^2\beta_j t)}e^{-\lambda_{j}t}dt\right]=
 \frac{\lambda_{j}}{\lambda_{j}+(1/2)\sigma^2_j q(\beta_j-q)},
 $$
if the denominator is positive, and  $f_{j}(q)=\infty$ otherwise. 

 Clearly, $f_{j}(q)<\infty$, 
if  $q\in [0,r_j[$,  $f_{j}(q)=\infty$, if $q\in [r_j,\infty[$, and $f_{j}(r_j-)=\infty$, 
where $r_j$ is the positive root of the equation 
 $$
 q^2-\beta_j q-2\lambda_{j}\sigma_j^{-2}=0,
 $$
 that is, $r_j=r(\lambda_j,\beta_j,\sigma_j)$, 
 \beq
 \label{r}
 r(\lambda,\beta,\sigma):=\beta/2+\sqrt{\beta^2/4+2\lambda \sigma^{-2}}. 
 \eeq
 

\smallskip 
Note that the formula (\ref{H0}) can be written in a shorter form  
\beq 
\label{H}
\Ups(q)=\E\left [\sum_{k=2}^\infty f_0(q)f_{\vartheta_1}(q),..., f_{\vartheta_{k-1}}(q)I_{\{\varpi \ge  k,\, \vartheta_k=0\}}\right]. 
\eeq
\smallskip

If $q\le \beta_*:=\min_j \beta_j$, then all $f_j(q)\le 1$ and $\Ups(q)$ is dominated the probability of  return of $\theta$ to the initial state (equal to unit), that is,  $[0,\beta_*]\subseteq {\rm dom}\, \Ups$. Also, $f_j(\beta_*/2)<1$ for all $j$ and, therefore, $\Ups(\beta_*/2)<1$.   
Since we assume that any state of $\theta$
 can be reached from any other state, ${\rm dom}\, \Ups\subseteq [0,r_*[$ where $r_*:=\min_j r_j$.


More precise information gives the following lemma. 
 
 \begin{lemma}  We have ${\rm dom}\, \Ups=[0,r_*[$  and 
 $\lim_{q\uparrow r^*}\Ups (q)=\infty$. 
  \end{lemma}
{\sl Proof.}  To explain the idea let us consider first the case $K=3$. Regrouping terms in the formula (\ref{H0}) according to 4 pairs ``exit from $0$ to $l$, return back from $m$" we get the representation
\bean
\Ups (q)&=& \big[ P_{0,1}f_0(q)P_{1,0}f_1(q)+P_{0,2}f_0(q) P_{2,0}f_2(q)  +P_{0,1}P_{1,2}P_{2,0}f_0(q)f_1(q)f_2(q)\\
&&+P_{0,2}P_{2,1}P_{1,0}f_0(q)f_2(q)f_1(q)\big] \sum_{k=0}^{\infty}(P_{1,2}f_1(q)P_{2,1}f_2(q))^k.
\eean
Note that if $P_{1,2}f_1(q)P_{2,1}f_2(q)<1$, then 
$$
 \sum_{k=0}^{\infty}(P_{1,2}f_1(q)P_{2,1}f_2(q))^k=\frac {1}{1-P_{1,2}f_1(q)P_{2,1}f_2(q)}, 
$$  
otherwise the above sum is equal to infinity.  Thus, $\Ups $ is a product of two continuous functions with values in $\bar \R_+$, hence, has the same property and the result follows. 

\medskip
For a model with an arbitrary $K$ we  get the continuity of $\Ups$ from the continuity result for more general functions. 

\smallskip
Let us consider a subset  $A\subset \{0,1,...,K-1\}$. 
For $i,k\notin A $  we denote by $\Gamma^A_{ik}$ the set of vectors $(i,i_1,i_2,\dots,i_m,k)$, 
$i_j\in A$, $j=1,\dots,m$, $m\in \N$. 
The elements of $\Gamma^A_{ik}$  are interpreted as cuts of  sample paths of the Markov chain entering to $A$ from the state $i$, evolving in $A$, and living $A$ to the state $k$.

Putting $h_{i,j}(q)=P_{i,j}f_i(q)$ with the natural convention $0\times \infty=0$  
we associate with elements of $\Gamma ^A_{ik}$ the continuous  functions 
$$
q\mapsto   h_{i,i_1}(q)
h_{i_1,i_2}(q)\dots h_{i_{m-1},i_m}(q)h_{i_m,k}(q)
$$
with values in $\bar \R_+$ and consider the sum of all these functions
$$
U^A_{ik}\colon q\mapsto \sum   h_{i,i_1}(q)
h_{i_1,i_2}(q)\dots h_{i_{m-1},i_m}(q)h_{i_m,k}(q). 
$$

Since $f_j<1$ on the interval $]0,\beta_*[$, also   $U^A_{ik}<1$ on this interval. 

We show by induction that  $U^A_{ik}\colon \R_+\to \bar \R_+$ is a continuous  function with $U^A_{ik}(0)\le 1$. Since $\Ups=U^{A\setminus\{0\}}_{00}$, this gives the assertion of the lemma. 

The idea of the proof consists in  representing  $U^A_{ik}$ as a sum of finite number of positive continuous functions using an appropriate  partition 
of  $\Gamma^A_{ik}$. 
Namely, for  $i_1\in A$ and $n\ge 0$ we define  the set   
$$
\Delta_{ik}^{i_1,n}:= \{i\}\times\Gamma^{A\setminus\{i_1\}}_{i_1,i_1}\times\dots \times\Gamma^{A\setminus\{i_1\}}_{i_1,i_1}\times  \Gamma^{A\setminus\{i_1\}}_{i_1,k},
$$  
composed by the vectors with the first component $i$, followed by $n \ge 0$ blocks formed by vectors from $\Gamma^{A\setminus\{i_1\}}_{i_1,i_1}$, and completed by  vectors from $\Gamma^{A\setminus\{i_1\}}_{i_1,k}$.  Clearly,  the countable  family   $\Delta_{ik}^{i_1,n}$, $i_1\in A$, $n\ge 0$, is a partition of 
$\Gamma^A_{ik}$ and  
\beq
\label{usum}
U^A_{ik}(q)=\sum_{i_1\in A}h_{i,i_1}(q) U^{A\setminus \{i_1\}}_{i_1k}(q) \sum_{n=0}^\infty\left [U^{A\setminus \{i_1\}}_{i_1i_1}(q)\right]^n. 
\eeq
\smallskip

The result is obvious when $A$ is a singleton, i.e.  $|A|=1$. 
Supposing that the assertion is already proven for the case where $|A|=K_1-1$ we consider the case  where $|A|=K_1$. 
By the induction hypothesis $ U^{A\setminus {i_1}}_{i_1m}\colon \R_+\to \bar \R_+ $ is a continuous  function for every  $i_1\in A$,  $m\notin A\setminus {i_1}$.   
The result follows from (\ref{usum}) and  the  formula for the geometric series.   \fdem

\medskip 
The strictly convex function $\Ups$ is less or equal to  unit on $[0,\beta_*]$, finite on $[0,r_*[$  and tends to infinity at $r_*$.  Hence, there is unique $\gamma\in ]\beta_*,r_*[$ such that 
 $\Ups(\gamma)=1$.  Moreover, $\Ups(\gamma+\epsilon)<\infty$ for some $\e>0$.  
 \smallskip

\smallskip

\section{Integrability of $Q_1$}

The following identity is obvious: 
$$
Y_{\ups_{n}}=-\sum_{k=1}^n \prod_{j=1}^{k-1}e^{-(V_{\ups_{j}}-V_{\ups_{j-1}})}\int_{]\ups_{k-1},\ups_{k}]}e^{-(V_s-V_{\ups_{k-1}})}dP_s.
$$
Using the abbreviations
$$
Q_k:=- \int_{]\ups_{k-1},\ups_{k}]}e^{-(V_s-V_{\ups_{k-1}})}dP_s, \qquad M_j:=e^{-(V_{\ups_{j}}-V_{\ups_{j-1}})} 
$$
we rewrite it in a more transparent form as 
\beq
\label{QM}
Y_{\ups_{n}}=Q_1+M_1Q_2+M_1M_2Q_3+...+M_1...M_{n-1}Q_n. 
\eeq

Note that the   random variables  $\ups_k-\ups_{k-1}$, that is, the  lengths of intervals between successive returns to the initial state, form and i.i.d. sequence. 
The random variables $(Q_k,M_k)$ have the same law and are  
 independent on the $\sigma$-algebra $\sigma\{(M_1,Q_1),..., (M_{k-1},Q_{k-1})\}$. 
 
First, we study integrability properties of $Y$. For this we need a general lemma 
involving parameters $\beta,\sigma>0$ and an independent on $W$  exponential random variable $\tau$ with parameter $\lambda>0$.   

\begin{lemma} 
\label{q1}
Let $0<q<r(\lambda,\beta,\sigma)$ where $r(\lambda,\beta,\sigma)$ is given by (\ref{r}). 
Then 
\bean
C(q,\lambda,\beta,\sigma):=\E\left [ \left (\int_{0}^\tau e^{-(\sigma W_s+(1/2)\sigma^2\beta s)}ds\right)^q\right]<\infty.
\eean 
\end{lemma}
\noindent{\sl Proof.} Put $W^{(\sigma\beta/2)}_s:=W_s+(1/2)\sigma\beta s$.  Take 
$\rho,\rho'>1$ such that $1/\rho+1/\rho'=1$ and $\rho q<r(\lambda,\beta,\sigma)$.
Dominating the integrant by its supremum and using the H\"older 
inequality we get that   
$$
C(q,\lambda,\beta,\sigma)\le \E\left [ \tau^q \sup_{s\le \tau} e^{-\sigma q W^{(\sigma\beta/2)}_s}\right]\le 
\left (\E\left [ \tau^{q\rho'}\right]\right)^{1/\rho'}\left(\E\left [ \sup_{s\le \tau} e^{-\sigma \rho q W^{(\sigma\beta/2)}_s}\right]\right )^{1/\rho}. 
$$
Since an exponential random variable has moments of  any order, the first multiplier in the right-hand side is finite. According to the formula (1.2.1) in Ch. 2 of the reference book \cite{BS},  
$$
\E\left [\sup_{s\le \tau} e^{-\sigma \rho q W^{(\sigma\beta/2)}_s} \right]= \E\left [ e^{-\sigma \rho q \inf_{s\le \tau}W^{(\sigma\beta/2)}_s}\right]= \frac{ r(\lambda,\beta,\sigma)}{ r(\lambda,\beta,\sigma)-\rho q}<\infty. 
$$
The lemma is proven. \fdem

\smallskip
Note that the condition $\Ups(q)<\infty$ holds only $q<r_*:=\min_jr(\lambda_j,\beta_j,\sigma_j)$ and the above lemma implies the following useful 
\begin{corollary} 
\label{coro}
If $\Ups(q)<\infty$, then $C^*(q):=\max_j C(q,\lambda_j,\beta_j,\sigma_j)<\infty$. 
\end{corollary} 

\begin{lemma} 
\label{L-int}
 Let $q>0$ be such that  $\Ups(q)<\infty$. Then 
\beq
\label{Leb}
\E\left [\int_0^{\ups_1} e^{-q V_s}ds\right]<\infty, \qquad 
\E \left [ \Big (\int_{0}^{\ups_1}e^{-V_s}ds\Big)^q\right]<\infty.
\eeq
\end{lemma}
{\sl Proof.}  Let $\tau$ be a random variable  exponentially distributed  with parameter 
$\lambda>0$ and let $W$ be a Wiener process independent on $\tau$.  Then 
\bean
\E\left [\int_{0}^\tau e^{-q (\sigma W_s+(1/2)\sigma^2\beta s)}ds\right] &=&
\int_{0}^\infty \P [s\le \tau] \E \left [e^{-q (\sigma W_s+(1/2)\sigma^2\beta s)}ds\right]\\
&=& \int_{0}^\infty e^{-(\lambda + (1/2)\sigma^2q (\beta - q))s}ds=
\frac {1}{\lambda + (1/2)\sigma^2q (\beta - q)},
\eean
if the denominator in the right-hand side  is strictly greater than  zero, and infinity otherwise. 

Using the   conditioning with respect to 
$\cF_{k-1}=\sigma \{\vartheta_1,...,\vartheta_{k-1}\}$ and the Markov property we get from (\ref{H}) with the abbreviation $\bar f_{k}(q):=f_0(q)f_{\vartheta_1}(q)... f_{\vartheta_{k-1}}(q)$, $k\ge 1$, that 
\bean
\Ups(q)&=& \sum_{k=2}^\infty \E \left [ \bar f_k(q)I_{\{\varpi \ge {k}, \vartheta_k=0\}}\right]=
\sum_{k=2}^\infty \E \left [ \bar f_k(q)I_{\{\varpi \ge {k}\}}\E [I_{\{ \vartheta_k=0\}}|\cF_{k-1}]\right]\\
&=&\E\left [ \sum_{k=2}^\infty  \bar f_k(q)I_{\{\varpi \ge {k}\}} \E [I_{\{ \vartheta_k=0\}}|\vartheta_{k-1}]\right ]\ge p_*\E\left [\sum_{k=2}^\infty  \bar f_k(q) I_{\{\varpi \ge {k}\}}\right] 
\eean 
where  $p_*>0$ is  the minimal of  values $P_{j,0}$ different from zero.  

Note that 
$$
\sum_{j=2}^{\varpi}\bar f_j(q)=\sum_{j=2}^\infty I_{\{\varpi = j\}}\sum_{k=2}^j\bar f_k(q)
= \sum_{k=2}^\infty \bar f_k(q)\sum_{j=k}^\infty I_{\{\varpi = j\}}=
\sum_{k=2}^\infty \bar f_k(q)   I_{\{\varpi \ge {k}\}}. 
$$
It follows  that 
$$
\E\left [ \sum_{j=2}^{\varpi}\bar f_j(q)\right]
=\E\left [\sum_{k=2}^\infty \bar f_k(q)   I_{\{\varpi \ge {k}\}}\right]\le \frac 1{p_*}\Ups(q).
$$
Since  
\bean
\E\left [\int_0^{\ups_1} e^{-q V_s}ds\right]&=&\E\left [ \sum_{k=2}^\infty I_{\{\varpi=k\}}\sum_{j=0}^{k-1}
 e^{-q V_{\tau_j}}
\int_{\tau_j}^{\tau_{j+1}}e^{-q (V_s-V_{\tau_j})}ds\right]\\ 
&=&\E \left [\sum_{k=2}^\infty I_{\{\varpi=k\}}\sum_{j=1}^{k}\bar f_{j}(q)  \frac {1}{\lambda_{\vartheta_{j-1}}}\right]\le \frac {1}{\lambda_{*}} \left (f_0(q)+\E\left [ \sum_{j=2}^{\varpi}\bar f_j(q)\right]\right), 
\eean
where $\lambda_*:=\min_i \lambda_i$,   
we obtain from the above inequalities that  
\beq
\label{xxxx}
\E\left [\int_0^{\ups_1} e^{-q V_s}ds\right]\le  \frac {f_0(q)}{\lambda_*}+\frac 1{\lambda_*p_*}\Ups (q)<\infty. 
\eeq
The first integrability property of (\ref{Leb}) is proven. 

\medskip

To prove the second property of (\ref{Leb}) we start with  the case
 $q\le 1$.  Using the elementary inequality $(\sum |x_i|)^q\le \sum |x_i|^q$ and Corollary \ref{coro}
 we get that  
 \bean
\E \left [ \left (\int_{0}^{\ups_1}e^{-V_s}ds\right)^q\right] &\le&  
\E\left [ \sum_{k=2}^\infty I_{\{\varpi=k\}}\sum_{j=0}^{k-1}
 e^{-q V_{\tau_j}}
\left(\int_{\tau_j}^{\tau_{j+1}}e^{-(V_s-V_{\tau_j})}ds\right )^q\right]\\ 
&= & \E \left [\sum_{k=2}^\infty I_{\{\varpi=k\}}\sum_{j=0}^{k-1}\bar f_k(q)\E\left [ \left (\int_{\tau_j}^{\tau_{j+1}}e^{-(V_s-V_{\tau_j})}ds\right )^q\right]\right]\\
&\le & C^*(q)\E\left [ \sum_{k=2}^\infty I_{\{\varpi=k\}}\sum_{j=0}^{k-1}\bar f_k(q)\right]\\
&\le& 
C^*(q)\left( 1+f_0(q)+\E\left [ \sum_{k=2}^\varpi \bar f_k(q)\right]\right)\\
&\le& C^*(q)\left( 1+f_0(q)+\frac 1{p_*} \Ups (q) \right)<\infty. 
\eean

\smallskip

Let  $q>1$. Due to the continuity of  $\Ups$ there exists $q'>q$ such that $\Ups(q')<\infty$. Applying, first,  the H\"older inequality with the conjugate exponents $q$ and $p:=q/(q-1)$ and then the   Young  inequality for products with the conjugate 
 exponents $q'/q$ and  $q'/(q'-q)$,  we get that 
\bean 
\E \left [ \left (\int_{0}^{\ups_1}e^{-V_s}ds\right)^q\right]&\le& 
\E \left [ \int_{0}^{\ups_1} \ups_1^{q-1} e^{-qV_s}ds\right]\\
& \le& (1-q/q')\E \left [\ups_1^{q_1}\right]+(q/q')
\E\left [\int_0^{\ups_1}e^{-q'V_s}ds\right], 
\eean
where $q_1:=(q-1)q'/(q'-q)$. 
The expectation of the integral in the right-hand side is finite in virtue of the first   inequality in (\ref{Leb}), applied with $q'$.  It remains to recall that the first return time of the finite state Markov process $\theta$ has moments of any order. 

For the reader convenience we give the proof of the above fact. Take  
an arbitrary  $m>1$ and let denote by $\varrho_j:=\tau_{j}-\tau_{j-1}$ the interjump times  of the Markov process $\theta$. Recall that the conditional distribution of the vector $(\varrho_1,\dots,\varrho_k)$ given 
$\vartheta_1=i_1,\dots,\vartheta_{k-1}=i_{k-1}$ is the same as the distribution 
of the vector $\tilde \varrho_1,\dots,\tilde \varrho_k)$ with independent components 
having, respectively, exponential distributions with the parameters $\lambda_0$,
$\lambda_{i_1}$ \dots, $\lambda_{i_{k-1}}$. Using   the  
H\"older inequality (now for the sum) and this fact we get that   
\bean
\E\left [ \ups_1^{m}\right]&=&\E \left [\left(\sum_{j=1}^{\varpi} \varrho_j\right)^{m}\right]\le 
\E\left [ \varpi^{m-1} \sum_{j=1}^\varpi \varrho_j^{m}\right]
\\
&=&
\E\left [\sum_{k\ge 2} \sum_{i_1\neq 0,i_2\neq 0,...,i_k=0}I_{\{\vartheta_1=i_1,\dots,
\vartheta_{k-1}=i_{k-1}, \vartheta_{k}=0
\}}k^{m-1}\sum_{j=1}^{k}\varrho_j^{m}\right]\\
&=&\Gamma (m+1) \E \left [\varpi^{m-1}  \left(\sum_{j=1}^{\varpi} \frac 1{\lambda_{\vartheta_{j-1}} }\right)^{m}\right]\le \Gamma (m+1) \lambda_*^{-m} \E \left [\varpi^{m}\right]
\eean
where $\Gamma$ is the Gamma function and $ \lambda_*:=\min_j \lambda_j$. 
It remains to make a reference to the fact 
that the first return time $\varpi$ for the Markov chain $\vartheta$ has moments of any order,  see, e.g. \cite{F}, Ch. XV, exer. 18-20.  \fdem

\smallskip The following lemma provides the required integrability property of $Q_1$.
\begin{lemma}  
\label{Q1}
Suppose that $\Pi(|x|^{\gamma}):=\int |x|^{\gamma}\Pi(dx)<\infty$. Then  $\E\left [ |Q_1|^{\gamma}\right]<\infty$.
\end{lemma}
{\sl Proof.}
{\bf Case where  $\gamma \le 1$}. The inequality  $(|x|+|y|)^\gamma\le |x|^\gamma+|y|^\gamma$ allows us to check  separately finiteness of moments of the integral of $e^{-V}$ with respect to 
the Lebesgue mesure (this is already done, see the second property in (\ref{Leb})) and the integral with respect to  the jump component 
of the process $P$. The latter integral is just a sum.  Since the  jump measure $\pi(dt,dx)$ has the compensator $\tilde \pi(dt,dx)=\Pi(dx)dt$, we have that   
\bean
\E\left [ |e^{-V}x  *\pi_{\ups_1}|^\gamma\right]
&\le& \E\left [ e^{-\gamma V}|x|^\gamma  *\pi_{\ups_1}\right]=\E \left [e^{-\gamma V}|x|^\gamma  *\tilde \pi_{\ups_1}\right]\\
&\le& \Pi(|x|^\gamma)\E\left [\int_0^{\ups_1} e^{-\gamma V_s}ds \right]<\infty
\eean
in virtue of the first property in (\ref{Leb}).

\smallskip 

{\bf Case where  $\gamma> 1$.} Now we shall split integrals using the elementary  inequality 
$$
(|x|+|y|)^\gamma\le 2^{\gamma-1}(|x|^\gamma+|y|^\gamma). 
$$
Because of the second property in (\ref{Leb}), we need to consider only the integral with respect to the jump component of $P$. 
Note that $e^{-V}|x|  *\tilde \pi_{\ups_1}<\infty$. Then 
$$
\E \left [\left(e^{-V}|x| *\pi_{\ups_1}\right )^\gamma\right]\le  2^{\gamma-1}\left (\E \left [\big|e^{-V}|x|  *(\pi-\tilde \pi)_{\ups_1}\big|^\gamma\right] + 
\E\left [\left( e^{-V}|x|  *\tilde \pi_{\ups_1}\right)^\gamma\right]\right).
$$
Due to the first property in (\ref{Leb})
$$
\E\left [\left ( e^{-V}|x|  *\tilde \pi_{\ups_1}\right)^\gamma\right] \le \big(\Pi(|x|) \big)^\gamma \E\left [  \left (\int_{0}^{\ups_1}e^{-V_s}ds\right)^\gamma\right]<\infty.
$$
Let $I_s:=e^{-V}|x|  *(\pi-\tilde \pi)_s$. According to the Novikov inequalities\footnote{
See \cite{Novikov} and a discussion in  \cite{KP2020}.}
 with $\alpha=1$  
the moment of the order $\gamma>1$ of the random variable $I^*_t:=\sup_{s\le t}|I_s|$ admits the bound  
\bean
\E\left [ I^{*\gamma}_{\ups_1} \right] &\le &C_{\gamma,1}\left(\E \left [(e^{-V}|x|*\tilde \pi_{\ups_1} )^\gamma\right] + \E\left [ e^{-\gamma V}|x|^\gamma *\tilde \pi_{\ups_1}\right]  \right)\\
&\le & C'_{\gamma,1}\E\left [ \left(\int_0^{\ups_1} e^{-V_s}ds\right)^\gamma\right]+ 
C'_{\gamma,1}\E\left [\int_0^{\ups_1} e^{-\gamma V_s}ds\right] 
\eean
where $C'_{\gamma,1}:=C_{\gamma,1}(\Pi(|x|))^\gamma<\infty$,  $C''_{\gamma,1}:=C_{\gamma,1}\Pi(|x|^\gamma)<\infty$ due to our assumption. 
The both integrals in the right-hand side,  as we proved, are finite.  \fdem

\section{Study of the process $Y$}
\begin{lemma} The process $Y$ has the following properties: 

$(i)$ $Y_t$ converges almost surely as $t\to \infty$ to a finite random variable 
$Y_\infty$.

$(ii)$ $Y_\infty=Q_1+M_1Y_{1,\infty}$ where $Y_{1,\infty}$ is a random variable independent on $(Q_1,M_1)$ and having the same law as $Y_{\infty}$. 

$(iii)$ $Y_{\infty}$ is unbounded from above. 
\end{lemma}
{\sl Proof.} $(i)$ Take $p\in ]0,\gamma\wedge 1 [$. Then $r:=\E\left [M_1^p\right]<1$ and, 
 Lemma \ref{Q1}, $\E \left[|Q_1|^p\right ]<\infty$. It follows that 
$\E\left [ |Y_{\ups_{n+1}}-Y_{\ups_{n}}|^p\right]=\E\left [ M_1^p\dots M_n^p Q_{n+1}^p\right]=r^n \E\left [ |Q_1|^p\right]$ and, therefore, 
$$
\E\left [ \left(\sum_{n\ge 0} |Y_{\ups_{n+1}}-Y_{\ups_{n}}|\right )^p\right]\le \sum_{n\ge 0}  \E \left [|Y_{\ups_{n+1}}-Y_{\ups_{n}}|^p\right]<\infty.
$$ 
Thus, $\sum_n |Y_{\ups_{n+1}}-Y_{\ups_{n}}|<\infty$ a.s. implying that $Y_{\ups_{n}}$ converges a.s. to some finite random variable we shall denote $Y_{\infty}$. 

Let $Y_t^*:=\sup_{s\le t}|Y_s|$.  Then  
$$
\E \left [Y^{*p}_{\ups_1}\right]
\le 
c^p  \E\left [\left(\int_0^{\ups_1}e^{-V_s}ds\right)^p\right] +(\Pi(|x|))^p \E\left [\left(\int_0^{\ups_1}e^{-pV_s}ds\right)^p\right]<\infty.
$$
Put 
$$
\Delta_n=\sup_{v\in [\ups_{2n},\ups_{n+1}]}\Big\vert\int_{\ups_{n}}^v e^{-V_s}dP_s\Big\vert.
$$
Then
$$
\E\left [ \Delta_n^p\right]= \E \left [M_1^p\dots M_n^p \sup_{v\in [\ups_{n},\ups_{n+1}]}\left\vert\int_{\ups_{n}}^v e^{-(V_{s}-V_{\ups_{n}})}dP_s\right\vert^p\right]\le r^n\E\left [ Y^{*p}_{\ups_1}\right]
$$
and, therefore, for any $\e>0$
$$
\sum_{n\ge 0}\P [ \Delta_n\ge \e]\le e^{-p}\sum_{n\ge 0}\E\left [ \Delta_n^p\right]<\infty. 
$$
By the Borel--Cantelli lemma for all $\omega$ except a null-set $\Delta_n(\omega)\le \e$ for all $n\ge n(\omega)$. This implies that $Y_t$ converges a.s.  to the same limit as the sequence $Y_{\ups_n}$. 
\smallskip

$(ii)$  Rewriting (\ref{QM}) in the form 
$$
Y_{\ups_{n}}=Q_1+M_1(Q_2+M_2Q_3+...+M_2...M_{n-1}Q_n)
$$
and observing that the sequence of random variables in the parentheses converges 
almost surely to a random variable with same law as $Y_{\infty}$ and independent on $(Q_1,M_1)$ we get the needed assertion. 

$(iii)$ In virtue of $(ii)$ it is sufficient to check that  the set 
$\{Q_1\ge N,\ M_1\le 1/N\}$ is non-null whatever is $N\ge 1$. 

Recall that 
$$
Q_1= -e^{-V}x*\pi_{\ups_1}-c\int_0^{\ups_1}e^{-V_s}ds, 
$$
where $
dV_s=\sigma_{\theta_s}dW_s+(1/2)\sigma^2_{\theta_s}\beta_{\theta_s}ds$. 

We consider several 
cases.
 
$1)$ $c<0$.   Using conditioning with respect to $\theta$ we may argue as $\theta$
would be deterministic, i.e. assuming that  $V$ is a  process with 
a deterministic switching of parameters and $\ups_1$ is just a number, say, $t>0$.   
On the set $\{T_1>\ups_1\}$  we have 
$Q_1=-c\int_0^{\ups_1} e^{-V_s}ds$. Since $T_1$ is independent on $W$ and the set $\{T_1>\ups_1\}$ we need to check only that 
the set 
$$
B_N(t):=\left\{-c\int_0^t e^{-V_s}ds\ge N,\ e^{-V_t}\le 1/N\right\}
$$ 
is non-null.  In the case where $\theta$ has no jumps on $[0,t]$ the process 
$V_s=\sigma_0W_t+(1/2)\sigma^2_{0}\beta_{0}t$ we get the latter property using conditioning with respect to $W_t=x$. Indeed,
 the conditional distribution of $(W_s)_{s\le t}$ given $W_t=x$ is the same as the (unconditional) distribution of the Brownian bridge $B^x=(B^x)_{s\le t}$ ending at $t$ at the value $x$. The latter is a continuous Gaussian process. This implies that 
the conditional distribution of the integral involved in $B_N(t)$ is unbound from above. Integrating over a suitable set with respect to the distribution of $W_t$ shows that  $B_N(t)$ is non-null.  In the case of several jumps at the moments $t_1,\dots t_k$ we can show that the integral over the interval $[0,t_1]$ has unbounded conditional distribution given $(W_{t_1},W_{t_2}-W_{t_1},\dots,W_t-W_{t_k})=(x_{t_1},x_{t_2},\dots x_{t_k+1})$ and conclude by integrating with respect 
to the distribution of the increments of $W$ over a set $[x,\infty[^{k+1}$ for sufficiently large $x\in \R_+$.   


$2)$ $c\ge 0$.   Put $\sigma^*:= \max_j \sigma_j $, 
$\kappa^*:=\max_j (1/2)\sigma_j^2\beta_j$, $\kappa_*:= \min_j (1/2)\sigma_j^2\beta_j$. Let $\delta>0$ and let 
 $r_N:=(\sigma^* \delta + \ln N)/\kappa_* $. The set  
$$
A_N:= \{ |W_s|<\delta,\ \forall\,s\le r_N+1\}\cap \{ r_N\le \ups_1\le r_N+1 \}
$$ 
 is non-null. On this set for all $s\in [0,\ups_1]$ we have the bounds
$$
- \sigma^* \delta+\kappa_* r_N\le V_s\le  \sigma^*\delta+\kappa^* (r_N+1)
$$
implying that 
$$
M_1=e^{-V_{\ups_1}}\le 
  e^{\sigma^*\delta-\kappa_* r_N} = 1/N
$$ 
and 
$$
c\int_0^{\ups_1}e^{-V_s}ds\le (c/N) (r_N+1)=:C_N. 
$$  
Since $P$ is not an increasing process, $\Pi(]-\infty,0[)>0$. Hence,  the set 
$$
\Big\{e^{-\sigma^*\delta - \kappa^* (r_N+1)}|x|I_{\{x< 0\}}*\pi_{r_N}\ge C_N +N,\ xI_{\{x >0\}}*\pi_{r_N}=0 \Big\}
$$
is non-null and its intersection with $A_N$ is also non-null. But this intersection  is a subset of the set $\{Q_1\ge N,\ M_1\le 1/N\}$. \fdem 

\section{Bounds for the ruin probability}
\label{bounds}
\begin{lemma}
\label{G-Paulsen} 
For every   $u>0$  
\beq
\label{Paulsen}
\bar G_i(u)\le\, \Psi_i(u)=\frac{\bar G_i(u)}{\E\big[\bar G_{\theta_{\tau^{u,i}}}(0) \vert \tau^{u,i}<\infty\big]}\le \frac{ \bar G_i(u)}{\min_{j}\bar G_j(0)},
\eeq  
where $\bar G_i(u):=\P[Y_\infty^i>u]$. 
\end{lemma}
\noindent
{\sl Proof.}
Let $\tau$ be an arbitrary stopping time with respect to the  filtration $(\cF^{P,R,\theta}_t)$. 
As the finite limit $Y^i_\infty$ exists,  the random variable 
$$
Y_{\tau,\infty}^i:=\begin{cases}
-\lim_{N\to \infty } 
\int_{]\tau,\tau+N]}\,e^{-(V_{s}-V_{\tau})}dP_{s},&\tau<\infty, \\
0, & \tau=\infty,
\end{cases} 
$$ 
is well defined.  On the set $\{\tau<\infty\}$
\beq
\label{YX}
Y_{\tau,\infty}^i=e^{V_\tau}(Y^i_{\infty}-Y^i_\tau)=X_{\tau}^u  +e^{V_\tau}(Y^i_\infty-u). 
\eeq   Let $\xi$ be a $\cF_{\tau}^{P,R,\theta}$-measurable random variable.  Note that  the conditional distribution of $Y_{\tau,\infty}^i$ given $(\tau,\xi,\theta_\tau)=(t,x,j)\in \bbr_+\times\bbr\times \{0,1,K-1\}$ is the same as the distribution of $Y_{\infty}^j$. It follows that   
$$
\P\big [
Y_{\tau,\infty}^i>\xi, \ 
\tau<\infty,\ \theta_\tau=j
\big] 
=\E\left [ \bar G_j(\xi)\, I_{\{ \tau<\infty,\; \theta_\tau=j\}}\right].
$$
Thus, if $\P[\tau<\infty]>0$, then
$$
\P\big[Y_{\tau,\infty}^i >\xi, \ 
\tau<\infty
\big]
=\E\big[\bar G_{\theta_\tau}(\xi)\, \vert\, \tau<\infty\big]\,\P\big[\tau<\infty\big]\,.
$$
Noting that  $\Psi_i(u):=\P\big[\tau^{u,i}<\infty\big]\ge \P\big[Y^i_{\infty}>u\big]=\bar G_i(u)>0$,   we deduce from here using (\ref{YX}) that 
\begin{align*}
\bar G_i(u)&=
\P\big[
Y^i_{\infty}>u,\ \tau^{u,i}<\infty\big]=
\P\big[Y_{\tau^{u,i},\infty}^i>X_{\tau^{u,i}}^{u,i},\ 
\tau^{u,i}<\infty
\big]\\
&=\E\big[\bar G_{\theta_{\tau^{u,i}}}(X_{\tau^{u,i}}^{u,i})\, \vert\, \tau^{u,i}<\infty\big]\,\P\big[\tau^{u,i}<\infty\big]
\end{align*}
implying the equality in (\ref{Paulsen}). Also, 
$$
\E\big[\bar G_{\theta_{\tau^{u,i}}}(0) \vert \tau^{u,i}<\infty\big]=
\sum_{j=0}^{K-1}\bar G_j(0)\P\big[\theta_{\tau^{u,i}}=j\vert \tau^{u,i}<\infty\big]\le {\min_{j}\bar G_j(0)} 
$$
implying the result. 
\fdem

\section{Ruin with probability one} 
\label{sec:prob1}

Assuming that $\beta^*:=\max_j\beta_j$ is strictly negative,  we give a sufficient condition  under which the ruin  is imminent.  

\begin{theorem}  
\label{Th.sec:prob1-1}
Suppose that $\beta^*<0$, 
 $\Pi(]-\infty,-\e])>0$ for all $\e>0$, and there exists  $\delta\in ]0,|\beta^*|\wedge 1[$ for which $\Pi(|x|^{\delta} )<\infty$.
Then $\Psi_i(u)=1$ for any $u>0$ and $i$.
  \end{theorem}
{\sl Proof.} Put $\wt{X}_{n}=\wt{X}^i_{n}:=X^i_{\upsilon_{n}}$. Note that  \eqref{uY} implies that the sequence $\wt{X}_{n}$ satisfies the  difference equation
\begin{equation}
\label{sec:prob1-1}
{\tilde X}_{n}={A}_{n}{\tilde X}_{n-1}
+
{B}_{n}, \qquad n\ge 1, \  {\tilde X}_{0}=u,  
\end{equation}
where ${A}_{n}:=M^{-1}_n:=e^{V_{\ups_{n}}-V_{\ups_{n-1}}}$ and 
$$
{B}_{n}:=-M^{-1}_n\,Q_{n}= \int_{]\ups_{n-1},\ups_{n}]}e^{V_{\ups_{n}}-V_s}dP_s.
$$
According to Corollary 6.2 in \cite{KP2020}, $\inf_n \wt{X}_{n}<0$ a.s.  if  the ratio $B_1/A_1$ is unbounded from below and  there is  $\delta\in ]0,1[ $ such that   $\E[A^{\delta}_{1}]<1$ and  $\E[|B_1|^{\delta}]<\infty$.  By our assumption  the event that 
on a fixed finite interval the process $P$ has arbitrary many   downward jumps of the size larger than $\e$ and has no jumps upward  is of strictly positive probability. Due to independence of $P$ and  $(W,\theta)$ this implies that $-Q_1=B_1/A_1$ is unbounded from below. 

Noting that 
\beq
\label{fltilde}
\tilde f_j(\delta):=\frac{\lambda_{j}}{\lambda_{j}+(1/2)\sigma^2_j \delta(|\beta_j|-\delta)}<1, 
\eeq
we get that 
$$
\E[A^{\delta}_{1}]=\E \big[e^{\delta V_{\ups_{1}}}\big]=\E\left [e^{\delta V_{\tau_1}}\prod_{i=2}^\varpi e^{\delta(V_{\tau_i}-V_{\tau_{i-1}})}\right]=\E\left[\sum_{i=2}^\varpi \tilde f_0(\delta )\tilde f_{\theta_1}(\delta )\dots  f_{\theta_{i-1}}(\delta )\right]<1. 
$$
Finally, the property  $\E[|B_1|^{\delta}<\infty$ can be proved by the same arguments as  in the  proof of 
Lemma \ref{Q1} with  $\gamma$ and $V$  replaced by $\delta$
and $V_{\ups_1}-V$ and the reference to (\ref{Leb}) replaced by the reference to (\ref{sec:prob1-IN-Q-2})
in Lemma \ref{L-int-22} below. \fdem 

\begin{lemma} 
\label{L-int-22}
Suppose that $\beta^*<0$. Then for any  $\delta\in ]0,|\beta^*|[$
\begin{equation}
\label{sec:prob1-IN-Q-2}
\E\left [\int_0^{\ups_1} e^{\delta(V_{\upsilon_1}- V_s})ds\right]<\infty, 
\quad\quad 
\E \left [ \left (\int_{0}^{\ups_1}e^{(V_{\upsilon_1}-V_s)}ds\right )^\delta\right]<\infty.
\end{equation}
\end{lemma}
{\sl Proof.}  The arguments are very similar to that of Lemma \ref{L-int} and we  only sketch them. The only new feature is that we need to consider processes of the form $(V_T-V_s)_{s\in [0,T]}$ rather than $(V_s)_{s\in [0,T]}$. 
The crucial observation is that the process $(W_T-W_s)_{s\in [0,T]}$ in the reversed time $s':=T-s$ is a Wiener process. 

First, observe that
$$
 \E\left [\int_0^{\ups_1} e^{\delta(V_{\upsilon_1}- V_s)}ds\right]=
\E\left[\sum^{\varpi}_{k=1} e^{\delta(V_{\upsilon_1}- V_{\tau_{k}})} \int_{\tau_{k-1}}^{\tau_{k}} e^{\delta(V_{\tau_{k}}- V_s)}ds \right]
$$
Given a trajectory of $\theta$,  the exponential and the integral in each summand are conditionally independent 
and their conditional expectations admit explicit expressions. For  the integral it is 
$1/\lambda_{\vartheta_{k-1}}\tilde f_{\vartheta_{k-1}}(\delta)$ where $\tilde f_j$  in (\ref{fltilde}).  
Note  that for $\delta\in ]0,|\beta^*|[$ the conditional 
expectation of the integral is dominated by $1/\lambda_*$ implying  that    
\bean
\E\left [\int_0^{\ups_1} e^{\delta(V_{\upsilon_1}- V_s)}ds\right]&\le &\frac 1{\lambda_*} \E\left[\sum^{\varpi}_{k=1} e^{\delta(V_{\tau_\varpi}- V_{\tau_{k}})}\right]\\
& =& \frac 1{\lambda_*} \E\left[\sum^{\varpi}_{k=1} \prod_{n=k}^{\varpi-1}e^{\delta(V_{\tau_{n+1}}- V_{\tau_{n}})}\right]
= \frac 1{\lambda_*}\E\left[\sum^{\varpi}_{k=1} \prod_{n=k}^{\varpi-1}\tilde f_{\theta_n}(\delta)\right] .
\eean
Due to the choice of $\beta$,  we have that $\tilde f^*:=\max_j \tilde f_j(\delta)<1$ and, therefore, 
$$
\E\left[\sum^{\varpi}_{k=1} \prod_{n=k}^{\varpi-1}\tilde f_{\theta_n}(\delta)\right]\le \E\left[\sum^{\varpi}_{k=1}
(\tilde f^*)^{\varpi-k}\right]\le \sum_{k=1}^\infty(\tilde f^*)^{k}<\infty. 
$$
The first property in (\ref{sec:prob1-IN-Q-2}) is proven.  

\smallskip

Let $\tau$ be an exponential random variable with parameter $\lambda>0$.  For any $\delta \in ]0,\wt r[$, where   
$$
\wt{r}:=\sqrt{{2\lambda}/{\sigma^2}+{\beta^{2}}/{4}}
+|\beta|/{2}, 
$$
we have, according to  (1.1.1) in Ch. 2 of the reference book \cite{BS}, that  
$$
\tilde C(\delta,\lambda,\beta,\sigma):=\E\left [ e^{\delta\sigma \sup_{s\le \tau}W^{(\sigma\beta/2)}_s} \right]= \frac{ \wt{r}}{\wt{r}-\delta}
<\infty.
$$
%
%
 We get (as in Corollary \ref{coro}) that for all $k\ge 1$
$$
\E\left[\left( \int_{\tau_{k-1}}^{\tau_{k}} e^{V_{\tau_{k}}- V_s}ds\right)^{\delta} \right]
\le \tilde C^*(\delta),
$$
with some constant $\tilde C^*(\delta)<\infty$, and we complete the proof of the second property in  (\ref{sec:prob1-IN-Q-2}) as in Lemma \ref{L-int}. \fdem 

\noindent

\smallskip
{\bf Acknowledgement.} This work was supported by the Russian Science Foundation associated grants  20-68-47030 and 20-61-47043.  


 \end{document}